\documentclass{article}

\usepackage{typearea}
\usepackage[latin1]{inputenc}
\usepackage{amsmath,amssymb,amsthm}
\usepackage[english]{babel}

\newcommand{\Z}{\ensuremath{\mathbb{Z}}}
\newcommand{\Zd}{\ensuremath{\mathbb{Z}^d}}

\newcommand{\R}{\ensuremath{\mathbb{R}}}

\newcommand{\E}{\ensuremath{\mathbb{E}\ }}

\newcommand{\Ld}{\ensuremath{\mathbb{L}^d}}

\newcommand{\Cov}{\text{Cov}}

\renewcommand{\epsilon}{\varepsilon}
\renewcommand{\limsup}{\overline{\lim}}
\renewcommand{\liminf}{\underline{\lim}}

\newcommand{\miniop}[3]{%
\renewcommand{\arraystretch}{0.6}
\begin{array}{c}
{\scriptstyle #1}\\
#2\\
{\scriptstyle #3}
\end{array}
\renewcommand{\arraystretch}{1}}

\newcommand{\Card}[1]{\vert #1 \vert}

\newcommand{\1}{1\hspace{-2.7mm}1}
\newcommand{\Var}{\text{Var }}

\newcommand{\supp}{\text{supp }}

\newtheorem{theo}{Theorem}
\newtheorem{lemme}{Lemma}
\newtheorem{coro}{Corollary}

\newtheorem{prop}{Proposition}

\author{Olivier \textsc{Garet}}

\date{Laboratoire de Mathématiques, Applications et Physique Mathématique
d'Orléans UMR 6628\\ Université d'Orléans\\ B.P. 6759\\
 45067 Orléans Cedex 2 France\\ E-Mail:
Olivier.Garet@univ-orleans.fr }

\title{Limit theorems for the painting of graphs by clusters\footnote{Preprint} }

\begin{document}

\maketitle

\centerline{\textbf{Abstract}}
%\foreignlanguage{english}
{We consider a generalization of the so-called \textbf{divide and
color model} recently introduced by Häggström
%as an alternative to
%Ising and Potts models
. We investigate  the behaviour of the
magnetization in large boxes and its fluctuations. Thus, laws of
large numbers and Central Limit theorems are proved, both quenched
and annealed. We show that the properties of the underlying
percolation process roughly influence the behaviour of the
colorying model. In the subcritical case, the limit magnetization
is deterministic and the Central Limit Theorem admits a Gaussian
limit. A contrario, the limit magnetization is not deterministic
in the supercritical case and the limit of the Central Limit
Theorem is not Gaussian, except in the particular model with
exactly two colors which are equally probable. }

\hspace{-5mm}AMS Classifications: 60K35, 82B20, 82B43.\\
KEY-WORDS: percolation, coloring model, law of large number,
central limit theorem.\\

\section{Introduction}
The aim of this paper is to give some results concerning a pretty
and natural model for the dependent coloring of vertices of a
graph. This model has been recently introduced by  Häggström
\cite{hag}, who presented the first results, concerning
essentially the presence (or absence) of percolation and the
quasilocality properties. The model is easily described: choose a
graph at random according to bond percolation, and then paint
randomly and independently the different clusters, each cluster
beeing monochromatic. There are several motivations for the study
of such a model, the most relevant being its links with Ising or
Potts models. We refer to the examples of the present article and
to the introduction of Häggström's paper for detailed motivations.

In Häggström's model, the panel was constituted by a finite number
of colours, which were chosen according to a measure $\nu$ on $\R$
with finite support. For our purpose, the natural assumptions will
only be the existence of a first or a second moment for $\nu$.

Actually, we will study the mean magnetization in large boxes: we
will identify the limit
$$M=\miniop{}{\lim}{n\to\infty}\frac1{\Card{\Lambda_n}}\sum_{x\in
\Lambda_n}X(x),$$
 and determinate its variations: we will prove central
 limit theorems for quantities such that
$$\frac1{(\Card{\Lambda_n})^{1/2}}\big( (\sum_{x\in \Lambda_n}
X(x))-\Card{\Lambda_n}M\big).$$ There are several natural
questions: when is $M$ deterministic ? What is the influence of
the underlying bond percolation ? When is there convergence to a
normal law in the  Central Limit Theorem ?

These questions can be asked in two different approach. We shall
use here the  vocabulary usually used in the theory of random
media.
\begin{itemize}
\item The quenched point of view: limit theorems are formulated
once the graph has been (randomly) fixed.
\item The annealed point of view: limit theorems are formulated
under the randomization of the graph.
\end{itemize}

Indeed, we will show that the properties of the underlying
percolation process roughly influence the behaviour of the
colorying model. In the subcritical case, the limit magnetization
is deterministic and the Central Limit Theorem admits a Gaussian
limit. A contrario, the limit magnetization is not deterministic
in the supercritical case and the limit of the Central Limit
Theorem is not Gaussian, except in the particular model with
exactly two colors which are equally probable. As examples, we
will study the case where $\nu$ is "+/-" valued and the case where
$\nu$ is a Gaussian measure.

\section{Notations}
We will deal here with stochastic processes indexed by $\Zd$.
Their definition will be related to some subgraphs of the
$d$-dimensional cubic lattice $\Ld$, which is defined by
$\Ld=(\Zd,E_d)$, where $E_d=\{ \{x,y\}\subset\Zd;
\sum_{i=1}^{d}\vert x_i-y_i\vert=1\}$.  In the following, the
expression "subgraph of $\Ld$" will always be employed for each
graph of the form $G=(\Zd,E)$ where $E$ is a subset of $E_d$. We
say that two vertices $x,y\in\Zd$ are adjacent in $G$ if
$\{x,y\}\in E$. Two vertices $x,y\in\Zd$ are said to be connected
in $G$ if one can find a sequence of vertices containing $x$ and
$y$ such that each element of the sequence is connected in $G$
with the next one. A subset $C$ of $\Zd$ is said to be connected
if each pair of vertices in $C$ are connected. The maximal
connected sets are called the connected components. They partition
$\Zd$. The connected component of $x$ and is denoted by $C(x)$.
Note that the connected components are also called clusters.
Conversely, a subset $D$ of $\Zd$ is said to be independent if no
pair in $D$ is constituted by adjacent vertices.

 We will consider here subgraphs of $\Ld$ which are generated
by Bernouilli bond percolation on $\Ld$. Thus, we will denote by
$\mu_p$ the image mesure of
$(\{0,1\}^{E_d},\mathcal{B}(\{0,1\}^{E_d}),((1-p)\delta_{0}+p\delta_{1})^{\otimes
E_d})$ by $$x\mapsto (\Zd,\{e\in E_d; x_e=1\}),$$ where $p\in
(0,1)$.

Let us choose a graph $G$ at random under $\mu_p$ and recall the
definition of some basic objects in percolation theory.

\begin{itemize}
\item The probability that $0$ belongs to an infinite cluster:\\ $\theta(p)=\mu_p(\Card{C(0)}=+\infty).$
\item The critical probability:\\ $p_c=\inf\{p\in (0,1); \theta(p)>0\}$.
\item The mean size of a finite cluster:\\ $\chi^f(p)=\sum_{k=1}^{+\infty} k
\mu_p(\Card{C(0)}=k).$
\item The number of open clusters per vertex\\ $\kappa(p)=\sum_{k=1}^{+\infty}
k^{-1} \mu_p(\Card{C(0)}=k).$
\end{itemize}

The following results will be currently used:
\begin{itemize}
\item $\mu_p$ is translation-invariant. As it is isomorphic to $(\{0,1\}^{E_d},\mathcal{B}(\{0,1\}^{E_d}),((1-p)\delta_{0}+p\delta_{1})^{\otimes
E_d})$, its tail $\sigma$-field is trivial and the ergodic theorem
can be employed with full power.
\item If $p\in (0,p_c)$, then $G$ contains no infinite cluster.
\item If $p\in (p_c,1)$, then $G$ contains $\mu_p$ almost surely one unique infinite
cluster.
\item If $p\ne p_c$, then $\chi^f(p)<+\infty$.
\end{itemize}

If $G$ is a subgraph of $\Ld$ and if $\nu$ is a probability
measure on $\R$, we will define the color-probability $P^{G,\nu}$
as follows: $P^{G,\nu}$ is the only measure on
$(\R^{\Zd},\mathcal{B}(\R^{\Zd}))$ under which the canonical
projections $X_i$ -- defined, as usually by $X_i(\omega)=\omega_i$
-- satisfy
\begin{itemize}
\item For each $i\in\Zd$, the law of $X_i$ is $\nu$.
\item For each independent set $S\subset\Zd$, the variables
$(X_i)_{i\in S}$ are independent.
\item For each connected set $S\subset\Zd$, the variables
$(X_i)_{i\in S}$ are identical.
\end{itemize}

The randomized color-measure is defined by $$P^{p,\nu}=\int
P^{G,\nu}\ d\mu_p(G).$$

 We also note
$\Lambda_n=\{-n,\dots,n\}^d$.

\section{Laws of large numbers}

\subsection{Quenched Law of large numbers}

\begin{theo}
\label{quenched_lln}
%For a given subgraph $G$ of $\Ld$, let us consider a family of
%random variables $(X^G(x))_{x\in\Zd}$ such that they are
%identically distributed and that $X(x)=X(y)$ if $x$ and $y$ are
%connected in $G$ and $X^G(x)$ is independent from $X^G(y)$ else.
%We suppose moreover that $X^G(0)$ is integrable and note $m=\E
%X^G(0)$.
Let $\nu$ be a probability measure on $\R$ with a first moment. We
put $m=\int_{\R} x\ d\nu(x)$. Let $p\in (0,1)\backslash\{p_c\}$

 For $\mu_p$ almost $G$, we have the following result %s
 :

%\begin{itemize}
%\item $G$ has at most one infinite connected component.
% \item
$$\miniop{}{\lim}{n\to
+\infty}\frac1{\Card{\Lambda_n}}\sum_{x\in
\Lambda_n}X(x)=(1-\theta(p))m+\theta(p)Z\quad P^{G,\nu}\text{
almost surely},$$ where $Z$ is the value taken by $X()$ along the
infinite component if it exists, and 0 else.
%\end{itemize}

\end{theo}
\begin{proof}

The following lemma will be of higher importance in the following.

\begin{lemme}
\label{quoca} For each subgraph $G$ of $\Ld$, let us denote by
$(A_i)_{i\in I}$ the partition of $G$ in connected component.

 Then, if $p\ne p_c$, we have
for $\mu_p$ almost $G$, we have $$\lim_{n\to\infty}\frac1{\Card
{\Lambda_n}}\sum_{i\in I;\Card{A_i}<+\infty}\Card{A_i\cap
\Lambda_n}^2=\chi^f(p),$$

where

$$\chi^f(p)=\sum_{k=1}^{+\infty} k P(C(0)=k).$$
\end{lemme}
\begin{proof}
Let us define $C'(x)$ by

\begin{equation*}
C'(x)= \begin{cases} C(x) & \text{ if }\Card{C(x)}<+\infty\\
\varnothing & \text{else}
\end{cases}
\end{equation*}
 and $C'_n(x)=C'(x)\cap \Lambda_n$.

It is easy to see that
\begin{equation}
\label{remarquable}
 \sum_{i\in I;\Card{A_i}<+\infty}\Card{A_i\cap
\Lambda_n}^2=\sum_{x\in \Lambda_n}\Card{C'_n(x)}.
\end{equation}
 We have
$C'_n(x)\le C(x)$, and the equality holds if and only if
$C'(x)\subset \Lambda_n$.

The quantity residing in connected components intersecting the
boundary of $\Lambda_n$ can be controlled using well-known results
about the distribution of the size of finite clusters. In both
subcritical case and supercritical case, we can found $K>0$ and
$\beta>0$ such that $$P(+\infty>\Card{C(x)}\ge n)\le\exp(-K
n^{\beta}).$$ (We can take $\beta=1$ when $p<p_c$ and
$\beta=(d-1)/d$ if $p>p_c$. See for example the reference book of
Grimmett \cite{MR2001a:60114} for a detailed historical
bibliography. ) It follows from a standard Borel-Cantelli argument
that for $\mu_p$ almost $G$, there exists a (random) $N$ such that
$$\forall n\ge N\quad \max_{x\in \Lambda_n} \Card{C'(x)}\le (\ln
n)^{2/\beta}.$$ If follows that for each $x\in \Lambda_{n-(\ln
n)^{2/\beta}}$, $C'(x)$ is completely inside $\Lambda_n$, and
therefore $C'(x)=C'_n(x).$ Then,
\begin{eqnarray*}
\sum_{x\in \Lambda_{n-(\ln n)^{2/\beta}}}\Card{C'(x)} & \le &
\sum_{x\in \Lambda_n}\Card{C'_n(x)}\le \sum_{x\in
\Lambda_n}\Card{C'(x)}\\ \frac1{\Card{\Lambda_n}}\sum_{x\in
\Lambda_{n-(\ln n)^{2/\beta}}}\Card{C'(x)} & \le &
\frac1{\Card{\Lambda_n}}\sum_{x\in \Lambda_n}\Card{C'_n(x)}\le
\frac1{\Card{\Lambda_n}}\sum_{x\in \Lambda_n}\Card{C'(x)}\\
\end{eqnarray*}

By the ergodic theorem, we have $\mu_p$ almost surely $$\lim_{n\to
+\infty} \frac1{\Card{\Lambda_n}}\sum_{x\in
\Lambda_n}\Card{C'(x)}=\E \Card{C'(0)}=\chi^f(p).$$ Since
$\lim_{n\to +\infty}\frac{\Card{\Lambda_{n-(\ln
n)^{2/\beta}}}}{\Card{\Lambda_n}}=1$ , the result follows.

\end{proof}

\noindent{\textbf{Remark}} If we forget technical controls, the
key point of this proof is the identity (\ref{remarquable}). It is
interesting to note that Grimmett \cite{MR2001a:60114} used an
analogous trick to prove that $\miniop{}{\lim}{n\to +\infty}
k(n)/\Card{\Lambda_n}=\kappa(p)$ almots surely, when $k(n)$ is the
number of open clusters in $\Lambda_n$.

Let $(a_i)_{i\ge 1}$ be a sequence such that for each $x\in\Zd$,
there exists exactly one $a_i$ connected to $x$ in $G$. Then
\begin{equation}
\label{decomposition} \sum_{x\in
\Lambda_n}X(x)=\miniop{+\infty}{\sum}{i=1}\Card{C'_n(a_i)}X(a_i)+Z\Card{\Lambda_n\cap
I},
\end{equation}
 where $I$ is the infinite connected component ($\varnothing$
if there is none). Since $$\Card{\Lambda_n\cap
I}=\miniop{}{\sum}{x\in \Lambda_n}\1_{\Card{C(x)}=+\infty},$$ it
follows from the ergodic theorem that
\begin{equation}
\label{in} \miniop{}{\lim}{n\to
+\infty}\frac1{\Card{\Lambda_n}}\Card{\Lambda_n\cap
I}=\E\1_{\Card{C(0)}=+\infty}=\theta(p) \end{equation}
 for $\mu_p$ almost $G$.
From now on, we will suppose that $G$ is a such a graph, and that,
moreover, it is such that the conclusions of lemma \ref{quoca}
hold -- $\mu_p$ almost all graph is such that.

Now, our goal is to apply a light improvement of the well-known
proof of Etemadi \cite{MR82b:60027} for the law of large number.
Let us state our result.

\begin{prop}
\label{etemadi} Let $(X_n)_{n\ge 1}$ be a sequence of pairwise
independent and identically distributed variables. We suppose that
$X_1$ is integrable and note $m=\E X_1$. Let also be
$(\alpha_{i,n})_{i,n\in\Z_{+}}$ be a doubly indexed sequence of
non-negative numbers such that
\begin{itemize}
\item $$\forall n\in\Z_{+}\quad k(n)=\Card{\{i\in\Z_{+};\alpha_{i,n}\ne 0\}}<+\infty.$$
%\item $$\lim_{n\to +\infty} k(n)=+\infty.$$
\item $$\forall i\in\Z_{+}\quad\text{the sequence } (\alpha_{i,n})_{n\ge
1}\quad\text{is a non-decreasing converging sequence}.$$
%\item $$\lim_{n\to\infty}
%{\miniop{+\infty}{\sum}{i=1}\alpha_{i,n}}=+\infty.$$
\item
 $$\frac{\sum_{i=1}^{+\infty}
\alpha_{i,n}^2}{\big(\sum_{i=1}^{+\infty}
\alpha_{i,n}\big)^2}=O(\frac1{k(n)}).$$
\item $$\exists A>0,d\ge 1\quad k(n)\sim
An^d.$$
\item $$\exists B>0,e>0\quad \sum_{i=1}^{+\infty}
\alpha_{i,n}\sim Bn^e.$$
\end{itemize}
Then, almost surely $$\lim_{n\to
+\infty}\frac{\miniop{+\infty}{\sum}{i=1}\alpha_{i,n}X_i}{\miniop{+\infty}{\sum}{i=1}\alpha_{i,n}}=m.$$
\end{prop}

\begin{proof}
By linearity, it is sufficient to prove the theorem for
nonnegative random variables. Assume then that $X_n\ge 0$.
Moreover, we can assume without loss that $\alpha_{i,n}=0$ for
$i>k(n)$. This can be done by permuting columns of the matrix
$(\alpha_{i,n})$. Let $d(n)=\sum_{i=1}^{+\infty} \alpha_{i,n}.$ We
define $S_n=\sum_{i=1}^{k(n)}\alpha_{i,n} X_i$ and $Q_n=S_n/d_n$.
We also consider the truncated variables $Y_i=X_i\1_{X_i\le k(i)}$
and associated sums and quotients:
$S^{*}_n=\sum_{i=1}^{k(n)}\alpha_{i,n} Y_i$ and
$Q^{*}_n=S^{*}_n/d_n$.

\begin{eqnarray*}
  \Var S_n^{*} & = & \sum_{i=1}^{k(n)}\alpha_{i,n}^2 \Var Y_i\\
  & \le & \sum_{i=1}^{k(n)}\alpha_{i,n}^2 \E Y_i^2\\
 & \le & \sum_{i=1}^{k(n)}\alpha_{i,n}^2 \E X_i^2\1_{X_i\le k(i)}\\
 & \le & \big(\sum_{i=1}^{k(n)}\alpha_{i,n}^2\big) \E X_1^2\1_{X_1\le k(n)}\\
\end{eqnarray*}
It follows that there exists $K>0$ such that $$\forall n\ge 1\quad
\Var Q_n^{*}\le K\frac1{k(n)} \E X_1^2\1_{X_1\le k(n)}.$$ Now, fix
$\beta>1$ and define  $u_n$ to be the integer which is the closest
to $\beta^n$. Then $$\sum_{n=1}^{+\infty}\Var Q_{u_n}\le  K\E
X^2_1\sum_{n=1}^{+\infty}\frac1{k(u_n)}\1_{X_1\le k(u_n)}$$ But
since $k(u_n)\sim A\beta^{nd}$, it is easy to prove that
$$\sum_{n=N}^{+\infty}\frac1{k(u_n)}=O(\frac1{k(u_n)}).$$
 Then, there exists $C>0$ such that
$$\forall N\ge 1\quad\sum_{n=N}^{+\infty}\frac1{k(u_n)}\le C
\frac1{k(u_n)}.$$ Now,
\begin{eqnarray*}
\sum_{n=1}^{+\infty}\Var Q_{u_n}^{*} & \le & K\E
X^2_1\sum_{n:k(u_n)\ge X_1 }\frac1{k(u_n)}\\ & \le & KC\E X^2_1
\frac1{k(\inf\{ n;k(u_n)\ge X_1\})}\\ & \le & KC\E X^2_1
\frac1{X_1}=KC\E X_1<+\infty\\
\end{eqnarray*}
It follows from Chebyshev's inequality and the first
Borel-Cantelli lemma that $$Q_{u_n}^{*}-\E Q_{u_n}^{*}\to 0\text{
a.s.}$$ By monotone convergence, $$\lim_{n\to +\infty} \E
Y_n=\lim_{n\to +\infty} \E X_1\1_{X_1\ge k(n)}=\E X_1.$$ Let $N$
be such that $\vert\E Y_k-\E X_1\vert\le\epsilon$ for $k>N$. Then,
for $n\ge N$ $$\vert\E Q_{u_n}-\E X_1\vert\le\frac{\E X_1}
{d_n}\big(\sum_{i=1}^N \lim_{k\to +\infty}
\alpha_{i,k}\big)+\epsilon.$$ Then, $\limsup_{n\to +\infty}\vert
\E Q_n^{*}-\E X_1\vert\le\epsilon$, and since $\epsilon$ is
arbitrary, $\lim \E{Q_n^{*}}=\E X_1$. It follows that
$$Q_{u_n}^{*}\to \E X_1\text{ a.s.}$$ Now, we go back to not
truncated variables. Since $A'=\sup_n {n/k(n)}<+\infty$, we have
$$\sum_{n=1}^{+\infty} P(X_n\ne Y_n)=\sum_{n=1}^{+\infty}
P(X_n>k(n) )=\sum_{n=1}^{+\infty} P(X_1>k(n)
)\le\sum_{n=1}^{+\infty} P(X_1>A'n)\le A'\E X_1<+\infty$$ It
follows that for almost all $\omega$, there exists $n(\omega)$
such that $X_k(\omega)=Y_k(\omega)$ for $k\ge n(\omega)$. Then,
for $n\ge n(\omega)$ $$\vert
Q_{n}(\omega)-Q^{*}_{n}(\omega)\vert\le\frac1
{d_n}\sum_{i=1}^{n(\omega)}\big( \lim_{k\to +\infty}
\alpha_{i,k}\big) X_i(\omega)).$$ It follows that $$Q_{u_n}\to \E
X_1\text{ a.s.}$$ If $u_n\le k\le u_{n+1}$, then since
$(S_n)_{n\ge 1}$ is non-decreasing, we have
$$\frac{d_{u_n}}{d_{u_{n+1}}}Q_{u_n}\le
Q_k\le\frac{d_{u_{n+1}}}{d_{u_n}}Q_{u_{n+1}}.$$ Since
$\frac{d_{u_{n+1}}}{d_{u_n}}\to \beta^e$, it follows that
$$\frac1{\beta^e}\E X_1\le\liminf_{k\to +\infty}Q_k
\le\limsup_{k\to +\infty}Q_k\le{\beta^e}\E X_1$$ Since this is
true for each $\beta>1$, we have proved that $$\lim_{k\to
+\infty}Q_k=\E X_1\text{ a.s.}$$
\end{proof}

% \ref{etemadi}.
Our goal is to apply this result to the sequence $(\alpha_{i,n})$
defined by $\alpha_{i,n}=\Card{C'_n(a_i)}$. Since the sequences
$(C'_n(x))_{n\ge 1}$ are non-decreasing, so are the sequences
$(\alpha_{i,n})_{n\ge 1}$. Indeed, we have $\miniop{}{\lim}{n\to
+\infty} \alpha_{i,n}=\Card{C(a_i)}<+\infty$. Moreover,
\begin{equation}
\label{comper}
{\miniop{+\infty}{\sum}{i=1}\alpha_{i,n}}=\Card{\Lambda_n\backslash
I}.
\end{equation}
 Since the $\alpha_{i,n }$'s are natural numbers, it follows
that $k(n)$ is finite. In fact, $k(n)$ is the number of finite
components which intersect $\Lambda_n$. Together with the ergodic
theorem, (\ref{comper}) gives
\begin{equation}
\label{hum}  {\miniop{+\infty}{\sum}{i=1}\alpha_{i,n}}\sim
 (1-\theta(p))\Card{\Lambda_n}\sim (1-\theta(p))2^d n^d.
\end{equation}
 %which implies that
 %${\miniop{+\infty}{\sum}{i=1}\alpha_{i,n}}$ tends to infinity
 %when $n$ does.
As already mentioned,  Grimmett has proved that
\begin{equation}
\label{kappa} \miniop{}{\lim}{n\to +\infty}
k(n)/\Card{\Lambda_n}=\kappa(p).
\end{equation}
%(see for example \cite{MR2001a:60114} for a nice proof)
Then, we have
\begin{equation}
\label{kappa2} k(n)\sim\kappa(p)2^d n^d
\end{equation}
 We must now prove that
  $$\frac{\sum_{i=1}^{+\infty}
\alpha_{i,n}^2}{\big(\sum_{i=1}^{+\infty} \alpha_{i,n}\big)^2}
k(n)$$ is bounded. But
 $$\frac{\sum_{i=1}^{+\infty}
\alpha_{i,n}^2}{\big(\sum_{i=1}^{+\infty} \alpha_{i,n}\big)^2}
k(n)\sim (1-\theta(p))^{-2} \frac{\sum_{i=1}^{+\infty}
\alpha_{i,n}^2}{\Card{\Lambda_n}}\frac{k(n)}{\Card{\Lambda_n}} $$
Using the conclusions of lemma \ref{quoca} and (\ref{kappa}), we
get $$\miniop{}{\lim}{n\to +\infty}\frac{\sum_{i=1}^{+\infty}
\alpha_{i,n}^2}{\big(\sum_{i=1}^{+\infty} \alpha_{i,n}\big)^2}
k(n)=\frac{\chi^f(p)}{(1-\theta(p))^{2}}\kappa(p) ,$$ which
completes the checking of the assumptions. It follows that
\begin{equation}
\label{morceau2} \miniop{}{\lim}{n\to
+\infty}\frac1{\Card{\Lambda_n\backslash I
}}\miniop{+\infty}{\sum}{i=1}\Card{C'_n(a_i)}X(a_i)=m\text{ a.s}.
\end{equation}
Since $\Card{\Lambda_n\backslash I }\sim
(1-\theta(p))\Card{\Lambda_n}$, it comes from
(\ref{decomposition}),(\ref{in}) and (\ref{morceau2}) that
$$\miniop{}{\lim}{n\to +\infty}\frac1{\Card{\Lambda_n}}\sum_{x\in
\Lambda_n}X(x)=(1-\theta(p))m+\theta(p)Z\quad P^{G,\nu}\text{
almost surely}.$$
\end{proof}

We will now formulate an easy, but important corollary.
\begin{coro}
\label{pas_trivial_1}
\begin{itemize}
\item $Z$  is measurable with respect to the tail
$\sigma$-field $$\mathcal{T}=\miniop{}{\cap}{n\ge
1}\mathcal{F}_{\Zd\backslash\Lambda_n},$$ where $\mathcal{F}_S$ is
the $\sigma$-field generated by the $(X_i)_{i\in S}$.
\item For $p>p_c$, $\mathcal{T}$ is not trivial under $P^{G,\nu}$
for $\mu_p$ almost every $G$ as soon as $\nu$ is not a Dirac
measure.
\end{itemize}
\end{coro}
\begin{proof}
The first point is a consequence of the formula given in
theorem~\ref{quenched_lln} and the second point is a consequence
of the first one, because $Z$ is non constant as soon as $\nu$ is
not a Dirac measure.
\end{proof}

The fact that $Z$ is $\mathcal{T}$-measurable will be important
for the formulation of annealed results, because the environment
is forgotten once we have randomized under $\mu_p$. Indeed,
whereas the infinite component can not always be recovered, the
value of $X()$ along this component can.

\subsection{Annealed Law of large numbers}

In this case, the annealed theorem is an easy consequence of the
quenched one.

\begin{coro}
Let $\nu$ be a probability measure on $\R$ with a first moment. We
put $m=\int_{\R} x\ d\nu(x)$.

Then, for $p\in (0,1)\backslash\{p_c\}$
 :

%\begin{itemize}
%\item $G$ has at most one infinite connected component.
% \item
$$\miniop{}{\lim}{n\to +\infty}\frac1{\Card{\Lambda_n}}\sum_{x\in
\Lambda_n}X(x)=(1-\theta(p))m+\theta(p)Z\quad P^{p,\nu}\text{
almost surely},$$ where $Z$ is the value taken by $X()$ along the
infinite component if it exists, and 0 else.
%\end{itemize}

\end{coro}
\begin{proof}
Let $C=\{\miniop{}{\lim}{n\to
+\infty}\frac1{\Card{\Lambda_n}}\sum_{x\in
\Lambda_n}X(x)=(1-\theta(p))m+\theta(p)Z\}$. We have
$$P^{p,\nu}(C)=\int P^{G,\nu}(C)\ d\mu_p(G)=\int 1\ d\mu_p(G)=1,$$
since $P^{G,\nu}(C)=1$ for $\mu_p$ almost $G$.
\end{proof}

Of course, the existence of the annealed limit is not an hard
result: since $P^{p,\nu}$ is invariant under the translations, the
ergodic theorem ensures that
\begin{equation}
\miniop{}{\lim}{n\to +\infty}\frac1{\Card{\Lambda_n}}\sum_{x\in
\Lambda_n}X(x)=\E [X_0\vert\mathcal{T}]\text{ almost surely}.
\end{equation}
Indeed, the annealed law of large number could be rephrased in
\begin{equation}\label{reformule}
\E [X_0\vert\mathcal{T}]=(1-\theta(p))m+\theta(p)Z
\end{equation}

Here is the annealed result analogous to corollary
\ref{pas_trivial_1}. It can be seen as a consequence of
(\ref{reformule}).

\begin{coro}
For $p>p_c$, $\mathcal{T}$ is not trivial under $P^{p,\nu}$  as
soon as $\nu$ is not a Dirac measure.
\end{coro}

\subsection{Examples}

\subsubsection{"+/-" valued spin system}

It is the simplest models that we can study: only two values are
taken: "+1" and "-1", with probability $\alpha$ and $1-\alpha$. In
the terminology of Häggström \cite{MR2001b:60118}, it is denoted
as the $r+s$-state fractional fuzzy Potts model at inverse
temperature $-\frac12\ln (1-p)$, with $r=\alpha$ and $s=1-\alpha$.
This name refers to the fact that the fuzzy Potts model can be
realized using random clusters by an analogous painting procedure.
For more details, see Häggström \cite{MR2001b:60118}.

We can also remark that for $\mu_p$ almost $G$, we have

$$P^{G,(1-\alpha)\delta_{-1}+\alpha\delta_{+1}}=\alpha\miniop{}{\lim}{\beta\to+\infty}\mathcal{I}^{+}_{G,\beta,h}
+(1-\alpha)\miniop{}{\lim}{\beta\to+\infty}\mathcal{I}^{-}_{G,\beta,h},$$
where
$\miniop{}{\lim}{\beta\to+\infty}\mathcal{I}^{+}_{G,\beta,h}$
(resp.
$\miniop{}{\lim}{\beta\to+\infty}\mathcal{I}^{-}_{G,\beta,h}$)is
the Ising Gibbs measure on $G$ at inverse temperature $\beta$ with
the external field $h=\frac12\ln(\alpha/(1-\alpha))$ which is
maximal (resp. minimal) for the stochastic domination. Thus,
$$P^{p,(1-\alpha)\delta_{-1}+\alpha\delta_{+1}}=
\miniop{}{\lim}{\beta\to+\infty}\int\alpha\mathcal{I}^{+}_{G,\beta,h}
+(1-\alpha)\mathcal{I}^{-}_{G,\beta,h}\ d\mu_p.$$ In this sense,
we can say that $P^{p,(1-\alpha)\delta_{-1}+\alpha\delta_{+1}}$
arises at the zero temperature limit of an Ising model on a
randomly diluted lattice. For precise definitions and results
relative to Ising ferromagnets on random subgraphs generated by
bond percolation, see Georgii~\cite{MR83g:82041} and also the
recent article of Häggström, Schonmann and Steif~\cite{MR1797305}.

 If we choose $\nu=(1-\alpha)\delta_{-1}+\alpha\delta_1$ with
$\alpha\in (0,1)$, it follows that the magnetization is
\begin{eqnarray}
\label{magnetisme} M=\miniop{}{\lim}{n\to
+\infty}\frac1{\Card{\Lambda_n}}\sum_{x\in \Lambda_n}X(x)=
 \begin{cases} 2\alpha(1-\theta(p))+2\theta(p)-1 & \text{ with probability }\alpha\\
 2\alpha(1-\theta(p))-1 & \text{ with probability }1-\alpha\\
\end{cases}
\end{eqnarray}
When $p\in (0,p_c)$, the magnetization is deterministic. Moreover,
it follows from (\ref{magnetisme}) that the signum of the
magnetization is deterministic if and only if
$$\max(\alpha,1-\alpha)(1-\theta(p))\ge\frac12.$$ Note that if
$\theta(p)\ge\frac12$ the signum of the magnetization can not be
deterministic.

Note that in the case $\alpha=\frac12$, the annealed law of the
magnetization has been identified by Häggström (\cite{hag}
Proposition 2.1) using a spin-flip argument.

\subsubsection{A quenched Gaussian system}

Here we choose $\nu=\mathcal{N}(0,1)$. For each $G$, $P^{G,\nu}$
is a Gaussian measure. Here, we have

\begin{eqnarray}
\label{magnetism} M=\miniop{}{\lim}{n\to
+\infty}\frac1{\Card{\Lambda_n}}\sum_{x\in
\Lambda_n}X(x)=\theta(p)Z.
\end{eqnarray}

In other words, $M$ is almost surely null when $p<p_c$ and $M\sim
\mathcal{N}(0,\theta(p)^2)$ when $p>p_c$.

We emphasize that these large numbers theorems are valid both
quenched and annealed. This will not more be so simple for Central
Limit theorems.

\section{Central Limit Theorems}

\subsection{Quenched Central Limit Theorem}

\begin{theo}
\label{quenched} Let $\nu$ be a probability measure on $\R$ with a
second moment. We put $m=\int_{\R} x\ d\nu(x)$ and
$\sigma^2=\int_{\R} (x-m)^2\ d\nu(x)$.

 For $\mu_p$ almost $G$, we have the following results:

\begin{itemize}
\item  \emph{The subcritical case}\\ If $p\in (0,p_c)$, then
$$\miniop{}{\lim}{n\to
+\infty}\frac1{\Card{\Lambda_n}^{1/2}}\big(\sum_{x\in \Lambda_n
}(X(x)-m)\big)=\mathcal{N}(0,\chi^f(p)\sigma^2).$$
 \item \emph{The supercritical case}\\ If $p\in (p_c,1)$, then$$\miniop{}{\lim}{n\to
+\infty}\frac1{\Card{\Lambda_n}^{1/2}}\big(\sum_{x\in
\Lambda_n\backslash I }(X(x)-m)\big)=\mathcal{N}(0,
\chi^f(p)\sigma^2)$$ where $I$ is the infinite component of $G$.
\end{itemize}

\end{theo}
For simplicity, we will give the proof in the supercritical case
-- which contains the proof of the subcritical case.
\begin{proof}
$$\sum_{x\in \Lambda_n\backslash I
}(X(x)-m)=\miniop{+\infty}{\sum}{i=1}\Card{C'_n(a_i)}(X(a_i)-m)$$
Then $$\frac1{\Card{\Lambda_n}^{1/2}}\sum_{x\in
\Lambda_n\backslash I
}(X(x)-m)=\big(\frac{s_n^2}{\Card{\Lambda_n}}\big)^{1/2}\frac1{s_n}\miniop{+\infty}{\sum}{i=1}\Card{C'_n(a_i)}(X(a_i)-m),$$
with $$s_n^2=\miniop{+\infty}{\sum}{i=1}\Card{C'_n(a_i)}^2.$$

By lemma \ref{quoca}, we have for $\mu^p$ almost $G$ $\lim_{n\to
+\infty}\frac{s_n^2}{\Card{\Lambda_n}}=\chi^f(p)$.

Now, we have just to prove
\begin{equation}
\frac1{s_n}\miniop{+\infty}{\sum}{i=1}\Card{C'_n(a_i)}(X(a_i)-m)\Longrightarrow\mathcal{N}(0,\sigma^2).
\end{equation}
Therefore, we will prove that for $\mu^p$ almost $G$, the sequence
$Y_{n,k}=\Card{C'_n(a_i)}(X(a_i)-m)$ satisfies the Lindeberg
condition. For each $\epsilon>0$, we have
\begin{eqnarray*}
\sum_{k=1}^{+\infty}\frac1{s_n^2}\int_{\vert Y_{n,k}\vert\ge
\epsilon s_n} Y_{n,k}^2\ d P^{G,\nu} & = &
\sum_{k=1}^{+\infty}\frac{\Card{C'_n(a_k)}^2}{s_n^2}\int_{
\Card{C'_n(a_k)}\vert x\vert\ge \epsilon s_n} (x-m)^2\ d\nu(x)\\ &
\le & \int_{\vert x\vert\ge \frac{\epsilon}{\eta_n}} (x-m)^2\
d\nu(x),
\end{eqnarray*}
with $\eta_n=\frac{\sup_{k\ge 1}\Card{C'_n(a_k)}}{s_n}$. Then, the
Lindeberg condition is fulfilled if $\lim \eta_n=0$. But we have
already seen that $s_n\sim (\chi^f(p)\Card{\Lambda_n})^{1/2}$,
whereas $\sup_{k\ge 1}\Card{C'_n(a_k)}=O((\ln n)^{2/\beta})$. This
concludes the proof.

\end{proof}

\subsection{Annealed Central Limit Theorem}
\begin{theo}
\label{annealed} Let $\nu$ be a probability measure on $\R$ with a
second moment. We put $m=\int_{\R} x\ d\nu(x)$ and
$\sigma^2=\int_{\R} (x-m)^2\ d\nu(x)$. Let $p\in (p_c,1)$. We
emphasize that $G$ is randomized under $\mu_p$.

\begin{itemize}
\item  \emph{The subcritical case}\\ If $p\in (0,p_c)$, then
$$\miniop{}{\lim}{n\to +\infty}\frac1
{\Card{\Lambda_n}^{1/2}}\big(\sum_{x\in
\Lambda_n}X(x)-m\Card{\Lambda_n})\big)=\gamma$$ where
$$\gamma=\mathcal{N}(0,\chi^f(p)\sigma^2)$$
 \item \emph{The supercritical case}\\ If $p\in (p_c,1)$, then

$$\miniop{}{\lim}{n\to +\infty}\frac1
{\Card{\Lambda_n}^{1/2}}\big(\sum_{x\in
\Lambda_n}X(x)-((1-\theta(p))m+\theta(p)Z)\Card{\Lambda_n})\big)=\gamma$$
where
%$$\gamma=\int_{\R}\mathcal{N}(0,\chi^f(p)\sigma^2+(z-m)^2\sigma_p^2)\
%d\nu(z)$$
$$\gamma\text{ is the image of }
\mathcal{N}(0,\chi^f(p)\sigma^2)\times\mathcal{N}(0,\sigma_p^2)\times\nu\text{
by }(x,y,z)\mapsto x+y(z-m),$$ with
$$\sigma_p^2=\miniop{}{\sum}{k\in\Zd}(P(0\in I\text{and }k\in
I)-\theta(p)^2),$$ where $I$ is the infinite component of $G$.
\end{itemize}

\end{theo}

In the subcritical case, the Annealed  Central Limit Theorem is a
simple consequence of the Quenched Central Limit Theorem.

In order to prove this result in the supercritical case, we will
need a Central Limit Theorem related to the variations of the size
of the intersection of the infinite cluster with large boxes.

\begin{prop}
\label{newman} Under $\mu_p$, we have $$\frac{\Card{\Lambda_n\cap
I}-\theta(p)\Card{\Lambda_n}}{{\Card{\Lambda_n}^{1/2}}}\Longrightarrow\mathcal{N}(0,\sigma_p^2),$$
where $I$ is the infinite component of $G$.
\end{prop}
\begin{proof}
$${\Card{\Lambda_n\cap
I(\omega)}-\theta(p)\Card{\Lambda_n}}=\sum_{k\in\Lambda_n}
f(T^k\omega),$$ where $T^k$ is the translation operator defined by
$T^k(\omega)=(\omega_{n+k})_{n\in\Zd}$ and
$f=\1_{\{\Card{C(0)}=+\infty\}}-\theta(p)$. Moreover $f$ is an
increasing function and $\mu^p$ satisfies the FKG inequalities.
Then, $(f(T^k \omega))_{k\in\Zd})$ is a stationary random field of
square integrable satisfying to the FKG inequalities. Therefore,
according to Newman \cite{MR81i:82070}, the Central Limit Theorem
is true if we prove that the quantity $$\sum_{k\in\Zd}
\Cov(f,f\circ T^k)$$ is finite. But $\Cov(f,f\circ
T^k)=\Cov(Y_0,Y_k),$ with $Y_k=\1_{\ \Card{C(k)}<+\infty \}}$.
Now, $$Y_k=\sum_{n=0}^{+\infty} F_{n,k},\text{with
}F_{n,k}=\1_{\{\Card{C(k)}=n\}}.$$ Hence
\begin{eqnarray*}
\Cov(Y_0,Y_k) & = & \sum_{n=0}^{+\infty} \sum_{p=0}^{+\infty}
\Cov(F_{n,0}, F_{p,k})\\ & =  &
\sum_{n=0}^{+\infty}\big(\Cov(F_{n,0}, F_{n,k})+2 \sum_{p=0}^{n-1}
\Cov(F_{n,0}, F_{p,k})\big)\\ & =  & \sum_{n>\Vert
k\Vert/2-1}^{+\infty}\big(\Cov(F_{n,0}, F_{n,k})+2
\sum_{p=0}^{n-1} \Cov(F_{n,0}, F_{p,k})\big)
\\
& =  & \sum_{n>\Vert k\Vert/2-1}^{+\infty}\Cov(F_{n,0}, F_{n,k}+2
\sum_{p=0}^{n-1}F_{p,k}),
\\
\end{eqnarray*}
because $F_{n,0}$ and $F_{p,k}$ are independent as soon as $\Vert
k\Vert\ge p+n+2$. Since $F_{n,0}\ge 0$ and $0\le F_{n,k}+2
\sum_{p=0}^{n-1}F_{p,k}\le 2$, we have $$\vert \Cov(F_{n,0},
F_{n,k}+2 \sum_{p=0}^{n-1}F_{p,k})\vert\le 2\E
F_{n,0}=2P(\Card{C(0)}=n).$$

Then,
%\begin{eqnarray*}\vert
$\Cov(Y_0,Y_k) \vert  \le
\sum_{n>\Vert k\Vert/2-1}^{+\infty}2P(\Card{C(0)}=n)%\\ & \le &
%2P(\Vert k\Vert/2-1<\Card{C(0)}<+\infty)
%\end{eqnarray*}
$ and $$\sum_{k\in\Zd} \Cov(f,f\circ T^k)\le
2\miniop{+\infty}{\sum}{n=1}\Card{\Lambda_{2(n+1)}}P(\Card{C(0)}=n).$$

Since Kesten and  Zhang \cite{MR91i:60278} have proved the
existence of $\eta(p)>0$ such that $$\forall n\in\Z_{+}\quad
P(\Card{C(0)}=n)\le\exp(-\eta(p)n^{(d-1)/d}),$$ it follows that
the series converges. Of course, a so sharp estimate is not
necessary for our purpose. Estimates derived from Chayes, Chayes
and Newman \cite{MR88i:60161}, and from Chayes, Chayes, Grimmett,
Kesten and Schonmann \cite{MR91i:60274} would have been
sufficient.
\end{proof}

\begin{proof}
Rearranging the terms of the sum, we easily obtain
$$\big(\sum_{x\in
\Lambda_n}X(x)-((1-\theta(p))m+\theta(p)Z)\Card{\Lambda_n})\big)=\sum_{x\in
\Lambda_n\backslash I
}(X(x)-m)+(Z-m)(\Card{I\cap\Lambda_n}-\Card{\Lambda_n}\alpha)$$ We
will now put $$Q_n=\frac1 {\Card{\Lambda_n}^{1/2}}\big(\sum_{x\in
\Lambda_n}X(x)-((1-\theta(p))m+\theta(p)Z)\Card{\Lambda_n})\big),$$
and define $$\forall t\in\R\quad \phi_n(t)=\E \exp(iQ_n)$$ and
$$\forall t\in\R\quad \phi_{n,z}(t)=\E \exp(iQ_n)\vert\{Z=z\}$$ As
usually, it means that $\E (\exp(iQ_t)\vert Z)=\phi_{n,z}(Z)$. It
is also important to emphasize that the following properties are
fulfilled under $P^{p,\nu}$:
\begin{itemize}
\item  $G$ is independent from $Z$.
\item  $(X_k \1_{k\notin I})_{k\in\Zd}$ is independent from $Z$.
\end{itemize}

Therefore, we have $$\phi_{n,z}(t)=\E \exp(-
\frac{it}{\Card{\Lambda_n}^{1/2}}\sum_{x\in \Lambda_n\backslash I
}(X(x)-m)+(z-m)(\Card{I\cap\Lambda_n}-\Card{\Lambda_n}\alpha))$$
Conditioning by $\sigma(G)$ and using the fact that $I$ is
$\sigma(G)$-measurable, we get $\phi_{n,z}(t)=\E f_n(t,.)
g_{n}((z-m)t,.)$, with
\begin{eqnarray*}
f_n(t,\omega)& = & \E \exp(-
\frac{it}{\Card{\Lambda_n}^{1/2}}\sum_{x\in \Lambda_n\backslash I
}(X(x)-m)\vert\sigma(G)\\ & = & \int \exp(-
\frac{it}{\Card{\Lambda_n}^{1/2}}\sum_{x\in \Lambda_n\backslash
I(\omega) }(X(x)-m)\ dP^{G(\omega),\nu} \end{eqnarray*}
 and
\begin{eqnarray*}
g_{n}(t,\omega)& = &  \exp(-
\frac{it}{\Card{\Lambda_n}^{1/2}}(\Card{I(\omega)\cap\Lambda_n}-\Card{\Lambda_n}\alpha))
\\ %& = & \int \exp(-
%\frac{it}{\Card{\Lambda_n}^{1/2}}(\Card{I(\zeta)\cap\Lambda_n}-\Card{\Lambda_n}\alpha))\
%d\mu_p(\zeta)
\end{eqnarray*}
By theorem~\ref{quenched} we have for each $t\in\R$ and
$P^{p,\nu}$ almost $\omega$: $\miniop{}{\lim}{n\to + \infty}
f_{n}(t,\omega) =\exp(-\frac{t^2}{2}\chi^f(p)\sigma^2)$ Then, by
dominated convergence $$\lim_{n\to +\infty}\E
(f_n(t,.)-\exp(-\frac{t^2}2\chi^f(p)\sigma^2))
g_{n}((z-m)t,.)=0.$$ Then
\begin{eqnarray*}
\lim_{n\to +\infty}\E f_n(t,.) g_{n}((z-m)t,.) & = & \lim_{n\to
+\infty}\exp(-\frac{t^2}2\chi^f(p)\sigma^2)\E g_{n}((z-m)t,.)\\ &
= &
\exp(-\frac{t^2}2\chi^f(p)\sigma^2)\exp(-\frac{t^2}2(z-m)^2\sigma^2_p)
\end{eqnarray*}
where the last equality follows from Proposition~\ref{newman}. We
have just proved that $$\miniop{}{\lim}{n\to\infty}
\phi_{n,z}(t)=\exp(-\frac{t^2}2(\chi^f(p)\sigma^2+(z-m)^2\sigma^2_p)).$$
Since $\phi_{n}(t)=\int \phi_{n,z}(t)\ d\nu(z),$ we get
\begin{eqnarray*}
\miniop{}{\lim}{n\to\infty} \phi_{n}(t) & = & \int
%\widehat{\mathcal{N}(0,\chi^f(p)\sigma^2+(z-m)^2\sigma^2_p)}(t)
\exp(-\frac{t^2}2(\chi^f(p)\sigma^2+(z-m)^2\sigma^2_p)) d\nu(z)\\
& = & \int \exp(itx) \ d\gamma(x).
\end{eqnarray*}
By the theorem of Levy , it follows that
$Q_n\Longrightarrow\gamma.$
\end{proof}

\subsection{Examples}

\subsubsection{"+/-" valued spin system}

If we choose $\nu=(1-\alpha)\delta_{-1}+\alpha\delta_1$ with
$\alpha\in (0,1)$, it follows that
\begin{itemize}
  \item In the subcritical case $p\in (0,p_c)$, then
  $$\gamma=\alpha\mathcal{N}(0,4\alpha(1-\alpha)\chi^f(p)).$$
  \item  In the supercritical case $p\in (p_c,1)$, then
$$\gamma=\alpha\mathcal{N}(0,4\alpha(1-\alpha)\chi^f(p)+4(1-\alpha)^2\sigma_p^2)+
(1-\alpha)\mathcal{N}(0,4\alpha(1-\alpha)\chi^f(p)+4\alpha^2\sigma_p^2).$$

\end{itemize}

\noindent\textbf{Remarks}\\
\begin{enumerate}
\item For the "+/-" valued spin system in the subcritical case,
the annealed Central Limit Theorem can be simply proved without
using the quenched one: since $\int \omega_k \ dP^{G,\nu}=m$ for
each $k$ and each $G$, it follows that the covariance of $X_0$ and
$X_k$ under $P^{p,\nu}$ is
\begin{eqnarray*}
\Cov (X_0,X_k)  & = & \int \big(\int (\omega_0-m)(\omega_k-m)\
dP^{G,\nu}\big)\ d\mu_p(G)\\ & = &\int\sigma^2\1_{\{k\in
C(0)\}}d\mu_p(G)\\ & = & \sigma^2 P(k\in C(0))
\end{eqnarray*}
Then,
\begin{eqnarray*}
\miniop{}{\sum}{k\in\Zd}\Cov (X_0,X_k)  & = &
\miniop{}{\sum}{k\in\Zd} \int\sigma^2\1_{\{k\in C(0)\}}d\mu_p(G)\\
& = & \sigma^2\int\miniop{}{\sum}{k\in\Zd}\1_{\{k\in
C(0)\}}d\mu_p(G)\\
\\ & = &
\sigma^2\int\Card{C(0)}d\mu_p(G)\\ & = & \sigma^2\chi(p),
\end{eqnarray*}
with $\chi(p)=\E\Card{C(0)}=\chi^f(p)+\theta(p)(+\infty).$
$\miniop{}{\sum}{k\in\Zd}\Cov (X_0,X_k)$ is a convergent series
when $p<p_c$ and a divergent one else.\\ In the subcritical case,
the theorem of Newman \cite{MR81i:82070} ensures that the Central
Limit is valid as soon as the translation-invariant measure
$P^{p,\nu}$ satisfy to the F.K.G. inequalities. Since Häggström
and Schramm \cite{hag} have proved the F.K.G. inequalities for the
"+/-" valued spin system, we get a simple proof for the annealed
Central Limit Theorem in this particular case.
\item
In the case were $\alpha=\frac12$,  $\gamma$ is a Gaussian measure
as well in the subcritical case
($\gamma=\mathcal{N}(0,\chi^f(p))$) as in the supercritical case
($\gamma=\mathcal{N}(0,\chi^f(p)+\sigma_p^2)$). It provides an
example where there is a classical Central Limit Theorem whereas
the "susceptibility" $\miniop{}{\sum}{k\in\Zd}\Cov (X_0,X_k)$ is
infinite.\\
 It is the "only" case with a Gaussian limit in the supercritical case,
as tell the following remark.

\item
If  $p\in (p_c,1)$, then $\gamma$ is Gaussian if and only if there
exist $a,b\in\R$ such that $\nu=\frac12(\delta_{a}+\delta_{b})$.
\begin{proof}
Using the characteristic function, it is easy to see that $\gamma$
is Gaussian if and only if  $\gamma'=\int
{\mathcal{N}(0,(z-m)^2\sigma^2_p)}(t)d\nu(z)$ does. Let us define,
for $k\in\Z_{+}$: $m_k=\int z^{2k}\ d\gamma(z)$  . We have
\begin{eqnarray*}
m_k & = & \int {\mathcal{N}(0,(z-m)^2\sigma^2_p)}(x\mapsto
x^{2k})d\nu(z)\\ & = & \int \frac{(2k)!}{k!2^k}(z-m)^{2k}
d\nu(z)\\ & = & \frac{(2k)!}{k!2^k}\int (z-m)^{2k} d\nu(z)
\end{eqnarray*}
By definition of $\gamma'$, $\gamma'$ is a symmetric measure. So
if $\gamma'$ if Gaussian, it is centered and we have $$\forall
k\in\Z_{+}\quad m_k= \frac{(2k)!}{k!2^k} m_1^k.$$ Then, we have
$$\forall k\in\Z_{+}\quad\int (z-m)^{2k} d\nu(z)= m_1^k.$$ If we
denote by $\nu'$ the image of $\nu$ by $z\mapsto (z-m)^{2}$, we
have $$\forall k\in\Z_{+}\quad\int_{\R_{+}} z^{k} d\nu'(z)=m_1^k$$
Then, we have $$\text{supp ess }\nu'=\lim_{k\to +
\infty}\big(\int_{\R_{+}} z^{k}
d\nu'(z)\big)^{1/k}=m_1=\int_{\R_{+}} z d\nu'(z)$$ It follow that
for $\nu'$ almost $z$, $z=\text{supp ess }\nu'$: $\nu'$ is a Dirac
measure. Therefore,
$\supp{\nu}\subset\{m-\sqrt{m_1},m+\sqrt{m_1}\}$. Since $m=\int z\
d\nu(z)$, we necessary have
$\nu(m-\sqrt{m_1})=\nu(m+\sqrt{m_1})=\frac12$ and then
$\nu=\frac12(\delta_{m-\sqrt{m_1}}+\delta_{m+\sqrt{m_1}})$.
\end{proof}

\end{enumerate}

\subsubsection{The quenched Gaussian system}

In the case $\nu=\mathcal{N}(0,1)$, theorem \ref{annealed} takes
the following form:

\begin{itemize}
\item  \emph{The subcritical case}\\ If $p\in (0,p_c)$, then
$$\miniop{}{\lim}{n\to +\infty}\frac1
{\Card{\Lambda_n}^{1/2}}\big(\sum_{x\in
\Lambda_n}X(x))\big)=\gamma$$ where
$$\gamma=\mathcal{N}(0,\chi^f(p)\sigma^2)$$
 \item \emph{The supercritical case}\\ If $p\in (p_c,1)$, then

$$\miniop{}{\lim}{n\to +\infty}\frac1
{\Card{\Lambda_n}^{1/2}}\big(\sum_{x\in
\Lambda_n}X(x)-\theta(p)Z\Card{\Lambda_n})\big)=\gamma$$ where
%$$\gamma=\int_{\R}\mathcal{N}(0,\chi^f(p)\sigma^2+(z-m)^2\sigma_p^2)\
%d\nu(z)$$
$$\gamma\text{ is the image of } \mathcal{N}(0,I_{\R^3})\text{ by
}(x,y,z)\mapsto (\chi^f(p))^{1/2}x+\sigma_p yz,$$ with
$$\sigma_p^2=\miniop{}{\sum}{k\in\Zd}(P(0\in I\text{and }k\in
I)-\theta(p)^2),$$ where $I$ is the infinite component of $G$.
\end{itemize}

In this case $\gamma$ is a Gaussian measure when $p<p_c$ whereas
$\gamma$ is a Gaussian chaos of order $2$ for $p>p_c$.

\bibliographystyle{alpha}
\nocite{*}
\bibliography{ccv}

\end{document}